\documentclass[12pt,a4paper]{article}
\usepackage{amsmath}
\usepackage{amsthm}
\usepackage{amssymb}
\usepackage{amsbsy}
\usepackage{amsfonts}
\usepackage{amstext}
\usepackage{amscd}
\usepackage[dvips]{epsfig}

\numberwithin{equation}{section}
\theoremstyle{plain}
\newtheorem{thm}{Theorem}[section]

\newtheorem{lem}[thm]{Lemma}
\theoremstyle{definition}

\newtheorem{rem}[thm]{Remark}
\newtheorem{defi}[thm]{Definition}

\begin{document}
\title{A simple proof for monotone CLT}
\author{Hayato Saigo \\ Graduate School of Science,  Kyoto University,\\  Kyoto 606-8502, Japan}
\date{}

\maketitle

\begin{abstract}
In the case of monotone independence, 
the transparent understanding of the mechanism to validate the central limit theorem (CLT) has been lacking,  
in sharp contrast to commutative, free and Boolean cases. 
We have succeeded in clarifying it by making use of simple combinatorial 
structure of peakless pair partitions.

\end{abstract}


\section{Introduction}
As a generalization of probability spaces, we define the notion of algebraic probability space. 
\begin{defi}(Algebraic probability space)
An algebraic probability space is a pair $(A,\varphi )$ consisting of a unital *-algebra $A$ and of a state $\varphi $ .  
\end{defi}
An element $a \in A$ is called an algebraic random variable.  For algebraic random variables, quantities such as 
$\varphi (a_1a_2...a_n)$ are called mixed moments.

In quantum probability, the notion of independence is 
understood as a reduction rule in calculation of mixed moments.  
 Among many different rules, commutative, free \cite{V1}, Boolean \cite{B, S-W} and monotone independence 
 \cite{Lu, Mur1, Mur2} are known as basic types \cite{B-S1, Mur0, Mur4, Fra1}.  
In the present paper, we focus on monotone independence, which is defined later.

We first introduce a key notion, 
a peak for a mapping (not necessarily homomorphism) between finite ordered set.

\begin{defi}\label{peak}
Let $S, T$ are finite ordered sets. 
An element $s \in S$ 
is said to be a peak of a map $f: S \longrightarrow T$ if  
$f(s) > f(s')$ holds for any $s'$ $(\neq s)$ which is 
next to $s$ (that is, there is no elements between $s$ and $s'$ with respect to the order on $S$).
\end{defi}

Now we introduce the notion of monotone independence. 

\begin{defi}(Monotone independence).
Let $\left\{A_{\lambda}; \lambda \in \Lambda  \right\}$ be a family of *-subalgebras of $A$, where the 
index set $\Lambda $ is equipped with a linear order $<$.

$\left\{A_{\lambda} \right\}$ is said to be monotone independent if
\[
\varphi (a_1...a_i...a_n)=\varphi (a_i)\varphi(a_1...a_{i-1}a_{i+1}...a_n)
  \]
holds 
for $a_i\in A_{\lambda _i}\backslash \textrm{\boldmath $C$}1$, whenever $i$ 
is a peak of the mapping $j \mapsto \lambda_j$. 
\end{defi}
 
For a monotone 
(or commutative, free, Boolean) independent family of *-subalgebras $\left\{A_{\lambda } \right\} $, 
the following impotant property holds, which is known as the singleton condition \cite{A-O, Obata} discussed by 
von Waldenfels in 1970s.

\begin{defi}(Singleton condition).
Let $(A, \varphi )$ be an algebraic probability space. If a family of *-subalgebras $\left\{A_{\lambda } \right\}$ 
satisfies the singleton condition if for any finite sequence $\lambda _1,...,\lambda _n \in \Lambda $ with 
a \textbf{singleton} $\lambda _s$, i.e. $\lambda _s\neq \lambda _i$ for all $i\neq s$,{}
\[\varphi (a_1...a_s...a_n)=\varphi (a_s)\varphi(a_1...a_{s-1}a_{s+1}...a_n) \]
holds when $ a_i \in A_{\lambda_i} $.
\end{defi}
   
This condition is essential for understanding asymptotic behavior of algebraic random variables. 
In fact, 
the following theorem holds \cite{A-O, Obata}.

\begin{thm}\label{fundamental}
Let $(A,\varphi)$ be an algebraic probability space and $(a_n)$ a sequence of elements of $A$ 
which satisfies the following: 

i) $a_i=a_i^{\ast }$,

ii) $\varphi (a_i)=0$,

iii) $\varphi (a_i^2)=1$,

iv) $\left\{a_i \right\}$ has uniform mixed moments, i.e., for $m\geq 1$,
\[
\sup \left\{\vert \varphi (a_{n_1}a_{n_2}\cdots a_{n_m})\vert \mid  n_1,n_2,\cdots ,n_m \in \textrm{\boldmath $N$} \right\} < \infty
\]

v) $\left\{ \textrm{\boldmath $C$} [a_i] \right\}$satisfies the singleton condition. 

Then, 
\[
M_m:=\lim _{N \rightarrow \infty} \varphi \left(\left(\frac{1}{\sqrt{N}}\sum^{\infty }_{n=1} a_n \right)^m\right)
\]
can be computed as follows:
\[
M_{2m+1}=0,\; \qquad M_{2m}=\lim _{N\rightarrow \infty }\frac{1}{N^m}\sum _{n\in \Pi(\left\{ 1,2,\cdots , 2m \right\},
\left\{ 1,2,\cdots , N\right\})}\varphi (a_{n_1}a_{n_2}\cdots a_{n_{2m}})
\]
Here, 
\[\Pi(S,T):=\left\{f:S \longrightarrow T \text{ such  that }  \left| f^{-1}(t)\right| 
=2 \text{ or }0 \right\}\]
for finite sets S and T.

\end{thm}

\section{A simple proof for the monotone CLT}
In the case of commutative, free and Boolean independence, the corresponding central limit theorems (CLTs) can be proved 
directly from Theorem \ref{fundamental} by using some simple combinatorics. 
The purpose of this section is to derive the mononotone CLT(Theorem \ref{the Monotone CLT} ) from Theorem 
\ref{fundamental} by using simple 
combinatorial argument.

To do this, we define a subset of $\Pi (S,T)$ for finite \textbf{ordered} sets $S$ and $T$, 
which we call the set of \textbf{peakless} pair partitions $\Pi_{0}(S,T)$ as follows: 

\[
\Pi_{0}(S,T):=\left\{f \in \Pi (S,T) \mid \text{ there is no 
peak of } f \right\}.
\]

\begin{rem}
For those who are familiar with the notion of ``monotone partition''
defined  by Muraki to classify quasi-universal products 
\cite{Mur0}, we note that a mapping $f:S \longrightarrow T$ 
which is an element of $\Pi (S,T)$  belongs to $\Pi_{0}(S,T)$ 
if and only if it is  a monotone partition, when $S,T$ are finite linearly ordered set. 
\end{rem}

It is easy to see that for the case of monotone independence, only 
peakless pair partions contribute and each contributions are equal to $1$, in the expression of moments in 
Thm.\ref{fundamental}. Hence the problem is reduced to counting the number of the elements
 in $\Pi_{0}(\left\{ 1,2,\cdots , 2m \right\},
\left\{ 1,2,\cdots , N\right\})$.

\begin{lem}\label{lemma}
\[\left|\Pi_{0}(\left\{1,2,\cdots , 2m \right\},
\left\{1,2,\cdots , N\right\})  \right| = \binom{N}{m}\times (2m-1)!! \] for $N\geq m\geq 1$.

\begin{proof}
First note that 
$f$ is an element of $\Pi_{0}(\left\{ 1,\cdots , 2m \right\},
\left\{1,\cdots , m\right\})$  
if and only if 
$f\mid_{\left\{1,\cdots , 2m \right\}\backslash \left\{i,i+1 \right\}}$ is an element of  
$\Pi_{0}(\left\{1,\cdots , 2m \right\}\backslash \left\{i,i+1 \right\},\left\{1,\cdots , m-1 \right\})$ 
for $i$ which satisfies $f^{-1}(m)=\left\{i, i+1 \right\}$. Then we have
\begin{eqnarray*}
\lefteqn{} \\
\left|\Pi_{0}(\left\{ 1,\cdots , 2m \right\},
\left\{ 1,\cdots , m\right\})  \right|
 =&\sum_{i} \left|\Pi_{0}(\left\{ 1,\cdots , 2m \right\}\backslash \left\{i,i+1 \right\},
\left\{ 1,\cdots , m-1\right\})  \right| & \\
 =&(2m-1)\times \left|\Pi_{0}(\left\{ 1,\cdots , 2(m-1) \right\},
\left\{ 1,\cdots , m-1\right\})  \right| & \\
\end{eqnarray*}
for $N\geq m \geq 2$, and hence,  $\left|\Pi_{0}(\left\{ 1,\cdots , 2m \right\},
\left\{ 1,\cdots , m\right\})  \right|=(2m-1)!!$ holds.

It is easy to see that
$
\left|\Pi_{0}(\left\{ 1,\cdots , 2m \right\},
\left\{ 1,\cdots , N\right\})  \right|=\binom{N}{m}\times \left|\Pi_{0}(\left\{ 1,\cdots , 2m \right\},
\left\{ 1,\cdots , m\right\})  \right|,
$ 
and we obtain the lemma.

\end{proof}
\end{lem}

The essence of the proof above can be understood as follows: 

First note that $\Pi(\left\{ 1,2,\cdots , 2m \right\},
\left\{ 1,2,\cdots , N\right\})$ represents all the possible ways to devide $2m$ balls into $m$ pairs and to paint 
each pairs by different $m$ colours chosen from $N$ colours.
Then, every peakless partition corresponding to the way how to paint balls can be generated uniquely 
by the following algorithm. Put $2m$ balls in line, choose $m$ colours from $N$ colours. Paint a couple of neighbouring 
balls by the highest colour. Then repeat this procedure, treating all the balls already painted and all the colours 
already used. 

Now it is easy to count the way of such colourings.  The number of the choices of $m$ colours from $N$ colours is 
$\binom{N}{m}$. The number of the choices of a couple of neighbouring balls painted by the highest colour is $2m-1$, 
and the 
number of the choices of a couple of neighbouring (ignoring the balls which have painted ) balls painted by the 
second highest colour is $(2m-2)-1=2m-3,\cdots ,$ and so on.

\begin{rem}
The proof of the lemma is the prototype of combinatorial argument in \cite{H-S}.
\end{rem}

Then we can prove the monotone CLT \cite{Lu, Mur1, Mur2, A-O, Obata}.

\begin{thm}\label{the Monotone CLT}
Let $(A,\varphi)$ be an algebraic probability space and $(a_n)$ a sequence of elements of $A$ 
which satisfies the following: 

i) $a_i=a_i^{\ast }$,

ii) $\varphi (a_i)=0$,

iii) $\varphi (a_i^2)=1$,

iv) $\left\{ \textrm{\boldmath $C$} [a_i] \right\}$  has uniform mixed moments, 

v) $\left\{ \textrm{\boldmath $C$} [a_i] \right\}$ is monotone independent. 

Then
\[
\lim _{N \rightarrow \infty} \varphi \left(\left(\frac{1}{\sqrt{N}}\sum^{\infty }_{n=1} a_n \right)^m\right)
=\frac{1}{\pi}\int _{-\sqrt{2}}^{\sqrt{2}}\frac{x^m}{\sqrt{2-x^2}}dx
\]
for $m=0,1,2,...$.
\begin{proof}
By the lemma \ref{lemma}, 
\[
M_{2m}=\lim _{N\rightarrow \infty }N^{-m} \binom{N}{m}(2m-1)!!=
\lim _{N\rightarrow \infty } N^{-m} \binom{N}{m} m! \times  \frac{(2m-1)!!}{m!}=\frac{(2m-1)!!}{m!}
\]
This is nothing but the $2m$-th moment of the standard arcsine law. 
\end{proof}

\end{thm}

\section*{Acknowledgements} 
The author was inspired by the series lectures given by Prof. Nobuaki Obata, to whom he is grateful 
for encouragements. He thanks Prof. Izumi Ojima and Mr. Ryo Harada for many useful comments to refine drafts.
He also 
thanks Prof. Shogo Tanimura, Mr. Takahiro Hasebe, 
Mr. Hiroshi Ando, and Mr. Kazuya Okamura for their interests and comments.


\begin{thebibliography}{122}
    \bibitem{A-O} L. Accardi and N. Obata, Foundations of Quantum Probability Theory (in Japanese), Makino Shoten Publ.(2003),  
    \bibitem{B} M. Bo\.{z}ejko, Uniformly bounded representations of free groups, J. Reine Angew. Math. \textbf{377} (1987), 170-186.
    \bibitem{Fra1} U. Franz, Unification of boolean, monotone, anti-monotone, and tensor independence and L\'{e}vy Processes, Math. Z. \textbf{243} (2003), 779-816. 
    \bibitem{B-S1} A. Ben Ghorbal and M. Sch\"{u}rmann, Non-commutative notions of stochastic independence, Math. Proc. Comb. Phil. Soc. \textbf{133} (2002), 531-561. 
    \bibitem{H-S} T. Hasebe and H. Saigo, The Monotone Cumulants, Arxiv:0907.4896 .  
    \bibitem{Lu} Y. G. Lu, An interacting free fock space and the arcsine law, Probability and Math. Stat. \textbf{17}, Fasc.1 (1997), 149-166. 
    \bibitem{Mur1} N. Muraki, Non-commutative Brownian motion in monotone Fock space, Commun. Math. Phys. \textbf{183} (1997), 557-570.
    \bibitem{Mur2} N. Muraki, Monotonic independence, monotonic central limit theorem and monotonic law of small numbers, Infin. Dimens. Anal. 
    Quantum Probab. Relat. Top. \textbf{4} (2001), 39-58. 
    \bibitem{Mur0} N. Muraki, The five independences as quasi-universal products, Infin. Dimens. Anal. Quantum Probab. Relat. Top. \textbf{5}, no. 1 (2002), 113-134.
    \bibitem{Mur4} N. Muraki, The five independences as natural products, Infin. Dimens. Anal. Quantum Probab. Relat. Top. \textbf{6}, no. 3 (2003), 337-371. 
    \bibitem{Obata} N. Obata, Notions of independence 
    in quantum probability and spectral analysis of graphs, Amer. Math. Soc. Transl. \textbf{223} (2008), 115-136.
    \bibitem{S-W} R. Speicher and R. Woroudi, Boolean convolution, in Free Probability Theory, ed. D. Voiculescu, 
Fields Inst. Commun., vol. 12 (Amer. Math. Soc., 1997), 267-280. 
    \bibitem{V1} D. Voiculescu, Symmetries of some reduced free product algebras, Operator algebras and their connections with topology and 
    ergodic theory, Lect. Notes in Math. \textbf{1132}, Springer (1985), 556-588. 
     
\end{thebibliography}
\end{document}